\documentclass[preprint,12pt]{elsarticle}

\usepackage[top=1in, bottom=1in, left=1in, right=1in]{geometry}
\usepackage{amssymb}
\usepackage{amsmath}
\usepackage{hyperref}
\hypersetup{hypertex=true,
	colorlinks=true,
	linkcolor=blue,
	anchorcolor=blue,
	citecolor=blue}
\usepackage{tabularx}
\usepackage{graphicx}
\usepackage{float}
\usepackage{booktabs} 
\usepackage{multirow} 
\usepackage{caption}
\usepackage{subcaption}
\DeclareCaptionLabelFormat{myformat}{\textbf{Fig.}~#2}
\DeclareCaptionLabelFormat{Table_fm}{\textbf{Table~#2}}

\journal{Applied Mathematics and Computation}

\begin{document}

\begin{frontmatter}


 \title{Two-relaxation-time regularized lattice Boltzmann model for convection-diffusion equation with variable coefficients}
 
\author{Yuan Yu\corref{cor1}}
\ead{yuyuan@xtu.edu.cn} 
\author{Zuojian Qin}
\author{Haizhuan Yuan}
\author{Shi Shu}
 
\cortext[cor1]{Corresponding author.}
\affiliation[aff0]{organization={School of Mathematics and Computational Science,
		Xiangtan University}, 
	city={Xiangtan}, 
	postcode={411105},
	country={China}}
\affiliation[aff1]{organization={National Center for Applied Mathematics in Hunan},
	city={Xiangtan}, 
	postcode={411105}, 
	country={China}}
\affiliation[aff2]{organization={Hunan Key Laboratory for Computation and Simulation in Science and Engineering, Xiangtan University},
             city={Xiangtan}, 
             postcode={411105},
             country={China}}

\begin{abstract}
In this paper, a new two-relaxation-time regularized (TRT-R) lattice Boltzmann (LB) model for convection-diffusion equation (CDE) with variable coefficients is proposed. Within this framework, we first derive a TRT-R collision operator by constructing a new regularized procedure through the high-order Hermite expansion of non-equilibrium. Then a first-order discrete-velocity form of discrete source term is introduced to improve the accuracy of the source term. Finally and most importantly, a new first-order space-derivative auxiliary term is proposed to recover the correct CDE with variable coefficients. To evaluate this model, we simulate a classic benchmark problem of the rotating Gaussian pulse. The results show that our model has better accuracy, stability and convergence than other popular LB models, especially in the case of a large time step.
\end{abstract}




\end{frontmatter}


\section{Introduction}\label{sec_intro}

Besides achieving massive success in the field of hydrodynamics\cite{aidun2010lattice,gabbana2020relativistic,kruger2017lattice}, the lattice Boltzmann method (LBM) has also shown great potential in solving nonlinear partial differential equations (PDEs)\cite{kruger2017lattice,chai2018lattice,doolen2019lattice,du2020lattice}. As a highly representative PDE, the convection-diffusion equation (CDE) is widely used to describe the transport of various scalars (such as mass, concentration, and energy) caused by the convection and diffusion processes in nature and industry\cite{bird2002transport}. Currently, the LBM has become a powerful tool for solving CDEs\cite{chai2013lattice,chai2020multiple,li2017lattice}, standing on an equal footing with finite difference\cite{mazumder2015numerical}, finite volume\cite{mazumder2015numerical}, and finite element methods\cite{roos2008robust}, not only due to its inherent programming simplicity, parallelism, and ability to handle complex boundaries with ease but also because it can generate exact d'Alembert solutions for the convection part, resulting in minimal numerical dissipation\cite{marie2009comparison}. 

The original lattice Boltzmann (LB) models for CDEs must hypothesize a constant convective velocity to exactly recover the macroscopic equations\cite{guo1999fully,van2000convection,chopard2009lattice}. But in many applications, such as thermal convection, solute transport, multiphase flow, and combustion, the convective velocity of the CDE should be a function of space and time rather than a constant. Whereas the velocity is non-constant, Chopard et al.\cite{chopard2009lattice} presented the precise form of the error in the macroscopic equations recovered. To eliminate this error while retaining the simple ``collision-propagation'' operation of LBM, a popular strategy is to add a time-derivative auxiliary term into the evolution equation\cite{shi2009lattice,wang2015regularized,zhao2020block,chai2020multiple,li2015lattice,chen2021multiple,shang2020discrete,wang2021modified,zhang2022hybrid,zhang2022advection}. However, as we explain later, these models have insufficient stability and accuracy when employing a large time step. The works of Chopard et al.\cite{chopard2009lattice} and Chai et al.\cite{chai2013lattice} have also indicated that these models with a space-derivative term have higher accuracy than models with a time-derivative term. Unfortunately, both of them need to introduce Navier-Stokes equations (NSEs) to achieve this substitution. To address this problem, a new first-order space-derivative auxiliary term instead of the previous time-derivative term is proposed in this paper.  

Moreover, the collision operator (CO) is also an important factor affecting the robustness of models. The simplicity and efficiency of the classic single-relaxation-time (SRT) Bhatnagar-Gross-Krook (BGK) CO made it very popular for a period of time, despite its stability and precision being widely criticized in handling convection-dominated problems\cite{perko2018single,yoshida2010multiple,chai2016multiple}. To confront the inherent limitations of the BGK CO, several different COs have been proposed within the scope of hydrodynamics, such as two-relaxation-time (TRT)\cite{ginzburg2008two}, raw-moment-based multiple-relaxation-time (MRT)\cite{lallemand2000theory}, central-moment-based MRT\cite{premnath2009incorporating},  entropic\cite{ansumali2002single}, and regularized\cite{latt2006lattice,coreixas2017recursive} COs. Although it is difficult to pick the best one, some widely agreed conclusions should be noted. Firstly, the TRT CO maintains a certain level of simplicity and effectiveness like the BGK CO, while also improving numerical stability and convergence like the MRT CO without introducing too many free parameters. Then, the BGK CO exhibits the viscosity-dependent numerical slip as enforcing the bounce-back (BB) boundary scheme, which diminishes the typical advantage of LBM that complex boundaries can be handled with a simple BB scheme. Fortunately, this notorious numerical slip can be eliminated by adopting the TRT model with a specific free relaxation parameter. In addition, the regularized lattice Boltzmann (RLB) model is also gaining more and more attention from researchers\cite{latt2006lattice,coreixas2017recursive,feng2019hybrid}. For example, the recursive RLB model was proposed for high Reynolds number (up to $10^6$) flow by Coreixas et al.\cite{coreixas2017recursive} and Mattila et al.\cite{mattila2017high}. Feng et al.\cite{feng2019hybrid} proposed a hybrid recursive regularized thermal lattice Boltzmann model for high subsonic compressible flows (Mach number up to 0.86). All of these studies indicate that RLB model is characterized by both excellent numerical stability and respectable efficiency\cite{ezzatneshan2019comparative,coreixas2020impact}.

When it comes to the RLB model for CDEs, some important works should be reviewed. Wang and his collaborators\cite{wang2015regularized,wang2021modified} first proposed and improved the RLB model for CDEs. Zhang et al.\cite{zhang2022hybrid}, inspired by the work of Jacob et al.\cite{jacob2018new}, introduced a hybrid scheme in the classic model to calculate the first-order non-equilibrium moment. These works have some pioneering significance, but their models are essentially SRT CO and thus suffer from the inherent limitations as mentioned earlier. To confront these limitations, a two-relaxation-time regularized (TRT-R) CO was derived in this paper through a unique approach unrelated to the work of Zhao and his coauthors\cite{zhao2020block}. 

Combining the first-order space-derivative auxiliary term with the TRT-R CO, we proposed a new model that is particularly suitable for solving CDE systems with variable coeffients. The Chapman-Enskog (CE) analysis shows that this model can accurately recover CDE without introducing the NSEs. Numerical experiments show that our model has better stability, accuracy, and convergence than the previous models. 

The rest of this paper is organized as follows. In Section \ref{section2}, Construction of the TRT-R LB model for the CDE with variable coefficients is presented. In Section \ref{section3}, numerical results of the TRT-R LB model are made. Finally, a brief conclusion is given in Section \ref{section4}. 
\section{The mathematical model}\label{section2}
\subsection{Governing equations}
The CDE with source term considered in this paper is given as follows:
\begin{gather}\label{eq_macro}
	\partial_t \phi+\partial_\alpha\left(\phi u_\alpha\right)=\partial_\alpha \left( \kappa \partial_\alpha \phi \right)+S,
\end{gather}
where $\phi$ is a scalar function of time $t$ and position $x_{\alpha}$ with the subscript $\alpha$ being the generic index, $\kappa$ is the diffusion coefficent, $u_{\alpha}$ is the non-constant convective velocity, and $S$ is the source term.

\subsection{Construction of the TRT-R CO}
To solve Eq.~(\ref{eq_macro}), the LB model proposed by Shi et al.\cite{shi2008new} is widely adopted, which possesses the following evolution equation given by
\begin{equation}\label{eq_ev_bgk}
	g_i\left(x_\alpha+e_{i \alpha} \Delta t, t +\Delta t\right)=
	\mathcal{C}_{i} +\Delta t G_i +\Delta t F_i+\frac{\Delta t^2}{2} \partial_t F_i,
\end{equation}
where $g_{i}$ is the probability distribution function, $\mathcal{C}_{i}$ is the CO, $G_i$ is the auxiliary term to guarantee the exact recovery of Eq.~(\ref{eq_macro}), and the first-order velocity moment of discrete source term
\begin{equation}\label{eq_Fi}
	F_i=w_i S+\left(1-\frac{1}{2 \tau_{p,1}}\right) w_i \frac{e_{i \alpha}}{c_s^2} u_\alpha S
\end{equation}
is used to introduce the source term $S$ more accurately. All the variables just mentioned are functions of time $t$ and space $x_{\alpha}$.  $\tau_{p,1}$ is the relaxation time associated with the diffusion coefficient $\kappa$, $w_{i}$ are the weight coefficients, $c_s$ is the lattice speed of sound, and $e_{i\alpha}$ is the discrete velocity in the $i$-th direction ranging from $0$ to $(q-1)$ for a given $d$-dimensional D$d$Q$q$ lattice model. The two-dimensional cases are considered for simplicity, and it is straightforward to extend to three dimensions, so the standard D2Q9 lattice model is adopted in this paper. For this lattice model,  $\{e_{i\alpha},i=0,...,8\}=\{(0,0),(\pm{c},0),(0,\pm{c}),(\pm{c},\pm{c})\}$, where the lattice speed $c=\frac{\Delta x}{\Delta t}$ with $\Delta x$ and $\Delta t$ being the spatial and temporal step respectively, $w_0=4/9$, $w_{1 \sim 4}=1/9$, $w_{5 \sim 8}=1/36$, and $c_s=c/\sqrt{3}$.

In the work of Shi et al.\cite{shi2008new}, the most popular Bhatnagar–Gross–Krook (BGK) collision operator was adopted and given by
\begin{equation}\label{eq_N_bgk}
	\mathcal{C}_{i,\text{BGK}}=g_i^{\text{eq}}+\left(1-\frac{1}{\tau_p}\right)   g_i^{\text{neq}},
\end{equation}
with the non-equilibrium distribution function being
\begin{equation}\label{eq_gneq}
	g_i^{\text{neq}}=g_i-g_i^{\text{eq}}, 
\end{equation}
where $g_{i}^{\text{eq}}$ is the equilibrium distribution function taking a third-order form \cite{coreixas2017recursive} as
\begin{equation}
	g_i^{\text{eq}}=w_i \phi\left\{\mathcal{H}_{i,0}+\frac{\mathcal{H}_{i, \alpha}}{c_s^2} u_\alpha+\frac{\mathcal{H}_{i, \alpha \beta}}{2 c_s^4} u_\alpha u_\beta+\frac{\mathcal{H}_{i, \alpha \beta \gamma}}{6 c_s^6} u_\alpha u_\beta u_\gamma\right\}.
\end{equation}

To improve the accuracy and stability of the LB model for nonlinear CDE, Wang et al.\cite{wang2015regularized} extended the popular regularized procedure from fluid dynamics to CDE, in which the regularized non-equilibrium
\begin{equation}\label{eq_Reg_gneq}
	\mathcal{R}(g_i^{\text{neq}})= w_i \frac{\mathcal{H}_{i, \alpha}}{c_s^2} \mathcal{B}_\alpha^{\mathrm{neq}},
\end{equation}
was used in place of the original non-equilibrium (Eq.~(\ref{eq_gneq})), where the first-order non-equilibrium moment $\mathcal{B}_\alpha^{\text{neq}}$ reads
\begin{equation}\label{eq_B1neq}
	\mathcal{B}_\alpha^{\mathrm {neq }}=\sum_i \mathcal{H}_{i, \alpha}g_i^{\text{neq}},
\end{equation}
with $\mathcal{H}_{i, \alpha}$ being the first-order Hermite polynomials tensor.
Thus, the Regularized CO becomes
\begin{equation}\label{eq8}
	\mathcal{C}_{i,\text{REG}}=g_i^{\text{eq}}+\left(1-\frac{1}{\tau_p}\right) \mathcal{R}(g_i^{\text{neq}}). 
\end{equation}

In order to further improve the stability and accuracy of the CO, we consider the following high-order Hermite expansion of non-equilibrium,
\begin{equation}\label{eq9}
	\mathcal{R}_{H}(g_i^{\text{neq}})\equiv w_i \mathcal{H}_{i,0} \mathcal{B}_{0}^{\mathrm{neq}}+w_i \frac{\mathcal{H}_{i, \alpha}}{c_s^2} \mathcal{B}_\alpha^{\mathrm{neq}}+w_i \frac{\mathcal{H}_{i, \alpha \beta}}{2 c_s^4} \mathcal{B}_{\alpha \beta}^{\mathrm{neq}}+w_i \frac{\mathcal{H}_{i, \alpha \beta \gamma}}{6 c_s^6} \mathcal{B}_{\alpha \beta \gamma}^{\mathrm{neq}}+\cdots,
\end{equation}
where the first four Hermite polynomials tensors are $\mathcal{H}_{i,0}=1$, $\mathcal{H}_{i, \alpha}=e_{i \alpha}$, $\mathcal{H}_{i, \alpha \beta}=e_{i \alpha} e_{i \beta}-c_s^2 \delta_{\alpha \beta}$, and $\mathcal{H}_{i, \alpha \beta \gamma}=e_{i \alpha} e_{i \beta} e_{i \gamma}-c_s^2\left(e_{i \alpha} \delta_{\beta \gamma}+e_{i \beta} \delta_{\gamma \alpha}+e_{i \gamma} \delta_{\alpha \beta}\right)$ respectively, with $\delta_{\alpha \beta}$ being the Kronecker delta, and the zeroth-order non-equilibrium moment satisfies $\mathcal{B}_{0}^{\mathrm{neq}}=\sum_i g_i^{\text{neq}}=0$ with the help of zeroth-order moment conservation relation. Furthermore, since the second and higher-order Hermite expansion terms in Eq.~(\ref{eq9}) have no effect on the recovery of CDE, which will be verified in the subsequent Chapman-Enskog (CE) analysis, we could propose a new multiple-relaxation-time regularized (MRT-R) collision operator as
\begin{equation}\label{eq10}
	\mathcal{C}_{i,\text{MRT-R}}=g_i^{\text{eq}}+\left(1-\frac{1}{\tau_{p,1}}\right) w_i \frac{\mathcal{H}_{i, \alpha}}{c_s^2} \mathcal{B}_\alpha^{\mathrm{neq}}+\left(1-\frac{1}{\tau_{p,2}}\right)w_i \frac{\mathcal{H}_{i, \alpha \beta}}{2 c_s^4} \mathcal{B}_{\alpha \beta}^{\mathrm{neq}}+\left(1-\frac{1}{\tau_{p,3}}\right)w_i \frac{\mathcal{H}_{i, \alpha \beta \gamma}}{6 c_s^6} \mathcal{B}_{\alpha \beta \gamma}^{\mathrm{neq}}+\cdots, 
\end{equation}
where $\mathcal{B}_{\alpha \beta}^{\text{neq}}$ and $\mathcal{B}_{\alpha \beta \gamma}^{\text{neq}}$ are the second- and third-order non-equilibrium moments respectively, and $\tau_{p,2}$ and $\tau_{p,3}$ are the free relaxation times. For simplicity, we can keep only one free relaxation time $\tau_{p,2}$ and thus obtain the TRT-R CO
\begin{equation}\label{eq11}
	\mathcal{C}_{i,\text{TRT-R}}=g_i^{\text{eq}}+\left(1-\frac{1}{\tau_{p,1}}\right) w_i \frac{\mathcal{H}_{i, \alpha}}{c_s^2} \mathcal{B}_\alpha^{\mathrm{neq}}+\left(1-\frac{1}{\tau_{p,2}}\right)w_i \frac{\mathcal{H}_{i, \alpha \beta}}{2 c_s^4} \mathcal{B}_{\alpha \beta}^{\mathrm{neq}},
\end{equation}
where the second-order non-equilibrium moment $\mathcal{B}_{\alpha \beta}^{\mathrm {neq }}$ is defined by
\begin{equation}\label{eq_B2neq}
	\mathcal{B}_{\alpha \beta}^{\mathrm {neq }}=\sum_i \mathcal{H}_{i, \alpha \beta}\left(g_i-g_i^{e q}\right).
\end{equation}
It needs to be stated here that the block triple-relaxation-time (B-TriRT) CO proposed by Zhao et al.\cite{zhao2020block} can be converted to the present TRT-R CO with $k_0=1$, $\mathbf{K}_1=\frac{1}{\tau_{p,1}} \mathbf{I}$, and $\mathbf{K}_2=\frac{1}{\tau_{p,2}} \hat{\mathbf{I}}$, where the symbols $k_0$, $\mathbf{K}$, $\mathbf{K}_2$, $\mathbf{I}$ and $\hat{\mathbf{I}}$ have been defined by Zhao et al.\cite{zhao2020block}. Therefore, theoretically speaking, the TRT-R CO should possess the ability to eliminate the numerical slip of the BB scheme.

Altogether, we provide a novel alternative approach to construct the MRT/TRT regularized procedure, and the TRT-R LB model for CDE possesses the following evolution equation
\begin{align}\label{eq_ev}
	g_i\left(x_\alpha+e_{i \alpha} \Delta t, t +\Delta t\right)=
	g_i^{e q}&+\left(1-\frac{1}{\tau_{p,1}}\right) w_i \frac{\mathcal{H}_{i, \alpha}}{c_s^2} \mathcal{B}_\alpha^{\mathrm{neq}} \notag \\
	&+\left(1-\frac{1}{\tau_{p,2}}\right) w_i \frac{\mathcal{H}_{i, \alpha \beta}}{2 c_s^4} \mathcal{B}_{\alpha \beta}^{\mathrm{neq}} +\Delta t G_i +\Delta t F_i+\frac{\Delta t^2}{2} \partial_t F_i,
\end{align}
where the discrete source term $F_{i}$ is given by Eq.~(\ref{eq_Fi}), and the auxiliary term $G_{i}$ will be presented below.
\subsection{The auxiliary term in evolution equation} \label{sec_auxiliary}
The form of auxiliary term determines whether the LB model can correctly recover the targeted equation. 

In the work of Shi et al.\cite{shi2008new}, they assumed that the convective velocity was constant, and the CDE could be recovered with $G_i=0$. When the non-constant convection velocity is considered, a popular auxiliary term is commonly adopted, which is denoted as
\begin{equation} \label{eq_t1}
	G_i=w_i \frac{e_{i \alpha}}{c_s^2}\left(1-\frac{1}{2 \tau_{p,1}}\right) \partial_t\left(\phi u_{\alpha}\right).
\end{equation}
This auxiliary term allows the LB model to recover the CDE accurately, but for some common scenarios, this calculation of the time-derivation can be effectively avoided.

In general, the convective velocity is often governed by an equation such as NSE, and the CDE must be solved coupled with the aforementioned equation in the same time step. In this case, the convection velocity $u_{\alpha}$ is usually viewed as a function of space $x_{\alpha}$ only in the same time step. Based on this assumption, we propose a new auxiliary term to enhance the model's robustness by circumventing the computing of the time-derivative term, which takes
\begin{equation}\label{eq_aux_new}
	G_i=w_i \frac{e_{i \alpha}}{c_s^2}\left(1-\frac{1}{2 \tau_{p,1}}\right) \phi u_\beta \partial_\beta u_\alpha.
\end{equation}
\subsection{Chapman-Enskog Analysis}
Here, we demonstrate that the LB model proposed above can accurately recover CDE by CE analysis technique.
Considering a small parameter $\varepsilon$, we introduce the following multiscale expansions:
\begin{equation}\label{eq_expansions}
	\begin{split}
		g_i=g_i^{(0)}+&\varepsilon g_i^{(1)}+\varepsilon^2 g_i^{(2)}, \quad 
		\partial_t=\varepsilon \partial_t^{(1)}+\varepsilon^2 \partial_t^{(2)}, \quad 
		\partial_\alpha=\varepsilon \partial_\alpha^{(1)}, \\
		&G_i=\varepsilon G_i^{(1)}, \quad 
		F_i=\varepsilon F_i^{(1)}, \quad 
		S=\varepsilon S^{(1)}.
	\end{split}
\end{equation}
With Eqs.~(\ref{eq_B1neq}), (\ref{eq_B2neq}), and (\ref{eq_expansions}), the non-equilibrium moments can be expanded in the following form
\begin{equation}\label{eq_expNEMs}
	\begin{split}
		\mathcal{B}_\alpha^{\mathrm {neq}}=\mathcal{B}_\alpha^{\mathrm {neq,(0)}}+\varepsilon \mathcal{B}_\alpha^{\mathrm {neq,(1)}}+\varepsilon^2 \mathcal{B}_\alpha^{\mathrm {neq,(2)}}, \\
		\mathcal{B}_{\alpha \beta}^{\mathrm{neq}}=\mathcal{B}_{\alpha \beta}^{\mathrm{neq,(0) }}+\varepsilon \mathcal{B}_{\alpha \beta}^{\mathrm {neq,(1)}}+\varepsilon^2 \mathcal{B}_{\alpha \beta}^{\mathrm{neq,(2)}}.
	\end{split}
\end{equation}
Applying the Taylor expansion to Eq.~(\ref{eq_ev}), we get
\begin{align}\label{eq_multi}
	g_i+\Delta t\left(\partial_t+e_{i \alpha} \partial_\alpha\right) g_i&+\frac{\Delta t^2}{2}\left(\partial_t+e_{i \alpha} \partial_\alpha\right)^2 g_i
	= g_i^{e q}+\left(1-\frac{1}{\tau_{p,1}}\right) w_i \frac{\mathcal{H}_{i, \alpha}}{c_s^2} \mathcal{B}_\alpha^{\mathrm{neq}} \notag \\
	&+\left(1-\frac{1}{\tau_{p,2}}\right) w_i \frac{\mathcal{H}_{i, \alpha \beta}}{2 c_s^4} \mathcal{B}_{\alpha \beta}^{\mathrm{neq}}+G_i \Delta t+\Delta t F_i +\frac{\Delta t^2}{2} \partial_t F_i.
\end{align}
Substituting Eqs. (\ref{eq_expansions}) and (\ref{eq_expNEMs}) into Eq.~(\ref{eq_multi}), one can get the $\varepsilon^0$, $\varepsilon^1$ and $\varepsilon^2$ scales' equations as follows:
\begin{subequations}
	\begin{gather}
		\label{eq10a}
		\begin{align}
			\varepsilon^{0}: \quad 
			g_i^{(0)}=g_i^{e q}+\left(1-\frac{1}{\tau_{p,1}}\right) w_i \frac{\mathcal{H}_{i, \alpha}}{c_s^2} \mathcal{B}_\alpha^{\mathrm{neq,(0)}}+\left(1-\frac{1}{\tau_{p,2}}\right) w_i \frac{\mathcal{H}_{i, \alpha \beta}}{2 c_s^4} \mathcal{B}_{\alpha \beta}^{\mathrm{neq,(0)}},
		\end{align}\\
		\label{eq10b}
		\begin{align}
			\varepsilon^{1}: \quad 
			g_i^{(1)}&+\Delta t \left(\partial_t^{(1)}+e_{i \alpha} \partial_\alpha^{(1)}\right) g_i^{(0)}=\left(1-\frac{1}{\tau_{p,1}}\right) \frac{\mathcal{B}_\alpha^{\mathrm{neq,(1)}}}{c_s^2} w_i \mathcal{H}_{i, \alpha}+\left(1-\frac{1}{\tau_{p,2}}\right) \frac{\mathcal{B}_{\alpha \beta}^{\mathrm{neq,(1)}}}{2 c_s^4} w_i \mathcal{H}_{i, \alpha \beta}\notag \\
			&+G_i^{(1)} \Delta t+F_i^{(1)} \Delta t,
		\end{align}\\
		\label{eq10c}
		\begin{align}
			\varepsilon^{2}: \quad 
			g_i^{(2)}&+\Delta t \partial_t^{(2)} g_i^{(0)} +\Delta t\left(\partial_t^{(1)}+e_{i \alpha} \partial_\alpha^{(1)}\right) g_i^{(1)}+\frac{\Delta t^2}{2}\left(\partial_t^{(1)}+e_{i \alpha} \partial_\alpha^{(1)}\right)^2 g_i^{(0)}\notag \\
			&=\left(1-\frac{1}{\tau_{p,1}}\right) w_i \frac{\mathcal{H}_{i, \alpha}}{c_s^2} \mathcal{B}_\alpha^{\mathrm{neq,(2)}}+\left(1-\frac{1}{\tau_{p,2}}\right) w_i \frac{\mathcal{H}_{i, \alpha \beta}}{2 c_s^4} \mathcal{B}_{\alpha \beta}^{\mathrm{neq,(2) }}+\frac{\Delta t^2}{2} \partial_t^{(1)} F_i^{(1)}.
		\end{align}
	\end{gather}
\end{subequations}
Taking the first- and second-order moments on both side of Eq.~(\ref{eq10a}), one can obtain 
\begin{gather}
	\mathcal{B}_\alpha^{\mathrm{neq,(0)}}=0, \quad 
	\mathcal{B}_{\alpha \beta}^{\mathrm{neq,(0)}}=0,
\end{gather}
which means
\begin{gather}
	g_i^{(0)}=g_i^{e q}.
\end{gather}
Summing Eq.~(\ref{eq10b}) over $i$, and multiplying Eq.~(\ref{eq10b}) by $e_{i\alpha}$ then summing it over $i$, we have 
\begin{subequations}
	\begin{gather}
		\label{eq12a}
		\begin{align}
			\partial_t^{(1)} \phi+\partial_\alpha^{(1)}\left(u_\alpha \phi\right)=S^{(1)},
		\end{align}\\
		\label{eq12b}
		\begin{align}
			\mathcal{B}_\alpha^{\mathrm{neq,(1)}}=-\Delta t \tau_{p,1} & \left\{\partial_t^{(1)}\left(\phi u_\alpha\right)+\partial_\beta^{(1)}\left(\phi u_\alpha u_\beta+\phi c_s^2 \delta_{\alpha \beta}\right)-\left(1-\frac{1}{2 \tau_{p,1}}\right) \phi u_\beta \partial_\beta^{(1)} u_\alpha\right. \notag \\
			& \left.-\left(1-\frac{1}{2 \tau_{p,1}}\right) u_\alpha S^{(1)}\right\}.
		\end{align}
	\end{gather}
\end{subequations}
Multiplying $\frac{\Delta t}{2}\left(\partial_t^{(1)}+e_{i \alpha} \partial_\alpha^{(1)}\right)$ on both sides of Eq.~(\ref{eq10b}), substituting the result into Eq.~(\ref{eq10c}) and then summing it over $i$, we can gain
\begin{align}
	\partial_t^{(2)} \sum_i g_i^{(0)}+\partial_\alpha^{(1)} 
	& \left(1-\frac{1}{2 \tau_{p,1}}\right) \mathcal{B}_\alpha^{\mathrm{neq,(1)}}+\frac{\Delta t}{2} \partial_\alpha^{(1)}\left\{\left(1-\frac{1}{2 \tau_{p,1}}\right) \phi u_\beta \partial_\beta^{(1)} u_\alpha\right\} \notag \\
	+ & \frac{\Delta t}{2} \partial_\alpha^{(1)}\left\{\left(1-\frac{1}{2 \tau_{p,1}}\right) u_\alpha S^{(1)}\right\}=0.
\end{align}
Combining the above equation with Eq.~(\ref{eq12b}), we can obtain
\begin{align}\label{eq23}
	\partial_t^{(2)} \phi-\partial_\alpha^{(1)}\left\{\left(\tau_{p,1}-\frac{1}{2}\right) \Delta t c_s^2 \partial_\alpha^{(1)} \phi\right\}=\Delta t \partial_\alpha^{(1)}\left(\tau_{p,1}-\frac{1}{2}\right)\left\{\partial_t^{(1)}\left(\phi u_\alpha\right)+u_\alpha \partial_\beta^{(1)}\left(\phi u_\beta\right)-u_\alpha S^{(1)}\right\}.
\end{align}
As mentioned in Section \ref{sec_auxiliary}, at each time step, the convective velocity is often considered as a time-independent quantity. So we assume $\partial_t^{(1)}u_\alpha=0$ here, and then combine it with Eqs.~(\ref{eq23}) and (\ref{eq12a}), and we have
\begin{gather}
	\label{Aeq15}
	\partial_t^{(2)} \phi-\partial_\alpha^{(1)}\left(\tau_{p,1}-\frac{1}{2}\right) c_s^2 \Delta t \partial_\alpha^{(1)} \phi=0.
\end{gather}
Taking Eq.~(\ref{Aeq15})$\times \varepsilon^2+$Eq.~(\ref{eq12a})$\times \varepsilon$, we get
\begin{gather}
	\partial_t \phi+\partial_\alpha\left(\phi u_\alpha\right)=\partial_\alpha \left( \kappa \partial_\alpha \phi \right)+S,
\end{gather}
where $\kappa=\left(\tau_{p,1}-\frac{1}{2}\right) c_s^2 \Delta t$ is the diffusion coefficient.

Notice that the value of the free relaxation parameter $\tau_{p,2}$ has no impact on the recovered equation, but it has a significant effect on the numerical stability. The optimal value of the free parameter is still an open question. Unless otherwise noted, this paper takes the ``magic'' parameter $\Lambda_p=\frac{1}{12}$ \cite{d2009viscosity} to regulate the free relaxation parameter $\tau_{p,2}$, where the magic parameter is defined as 
\begin{gather}
	\Lambda_p=(\tau_{p,1}-1/2)(\tau_{p,2}-1/2).
\end{gather}

\section{Numerical results}\label{section3}
To validate the proposed LB model, a classic benchmark problem of the convection and diffusion of a Gaussian pluse in a two-dimensional rigid-body rotating flow is considered\cite{zhou2016general,ding2020semi}. The governing equation and the initial condition are  
\begin{gather}
	\partial_t{\phi}-\partial_x {(y \phi)}+\partial_y{(x \phi)}=\partial_\alpha \left( \kappa \partial_\alpha \phi \right)+S, \quad x, y \in[-2 \pi, 2 \pi]
\end{gather}
and
\begin{gather}
	\phi_{0}\left(x,y\right)=\exp\left[-\left(x^2+3y^2\right)\right]
\end{gather}
respectively, where the diffusion coefficient $\kappa$ is a constant and will be specified later, and the source term is defined as
\begin{gather}
	S\left(x,y,t\right)=\left[6\kappa-4xy-4\kappa\left(x^2+9y^2\right)\right]\exp\left[-\left(x^2+3y^2+2\kappa t\right)\right].
\end{gather}
The corresponding analytical solution can be represented by\cite{ding2020semi}
\begin{gather} \label{eq_ana}
	\phi\left(x,y,t\right)=\exp\left[-\left(x^2+3y^2+2\kappa t\right)\right].
\end{gather}

To evaluate the accuracy, we define the root mean square (RMS) error as
\begin{equation}\label{eq19}
	E_{r}=\sqrt{\frac{\sum_{x,y}{\left|\phi_{\text{numeric}}-\phi_{\text{analytic}}\right|^2}}{N_x \times N_y} },
\end{equation}
where $N_x$ and $N_y$ denote the numbers of mesh in the $x$ and $y$ directions, respectively. In simulations, we set $N_x=N_y=N=L/\Delta x$ for the square-shaped domain with the characteristic length $L=4\pi$. $g_{i}=g_{i}^{\text{eq}}(\phi_{0},u_{\alpha})$ is applied to implement the initial condition, and the non-equilibrium extrapolation scheme is adopted to execute the boundary condition according to Eq.~(\ref{eq_ana})\cite{wang2021modified,wang2015regularized}. All cases were run until two rotations were completed. The RMS error at the end of the first rotation ($t=1$) was recorded, and we consider the case to be stable if it is still no blow-up upon the end of the second rotation ($t=2$). Unless otherwise specified, two other models were employed to simulate each example for comparison. In the first model, we specify $k_0=1$, $\mathbf{K}_1=\frac{1}{\tau_{p,1}} \mathbf{I}$, and $\mathbf{K}_2=\frac{1}{\tau_{p,2}} \hat{\mathbf{I}}$ in the model of Zhao et al.\cite{zhao2020block} to exclude the influence of CO and verify the advantages of the newly proposed auxiliary term and discrete source term (ADST) given by Eqs.~(\ref{eq_aux_new}) and (\ref{eq_Fi}) respectively. Hence, the difference between these two models is the ADST. In Zhao et al.'s model\cite{zhao2020block}, the conventional auxiliary term 
\begin{equation} \label{eq_t1}
	G_i=w_i \frac{e_{i \alpha}}{c_s^2}\left(1-\frac{1}{2 \tau_{p,1}}\right) \partial_t\left(\phi u_{\alpha}\right)
\end{equation}
and the traditional discrete source term
\begin{equation} \label{eq_t2}
	F_i=w_i S
\end{equation}
are adopted. Besides, the popular models\cite{chopard2009lattice,malaspinas2010lattice,su2013lattice} with the BGK CO and the traditional ADST (as given by Eqs.~(\ref{eq_t1}) and (\ref{eq_t2})) are also included for comparison.
\subsection{Accuracy and stability}
\begin{figure}[H]
 	\centering
 	\begin{subfigure}[b]{0.45\linewidth}
 		\includegraphics[width=\linewidth]{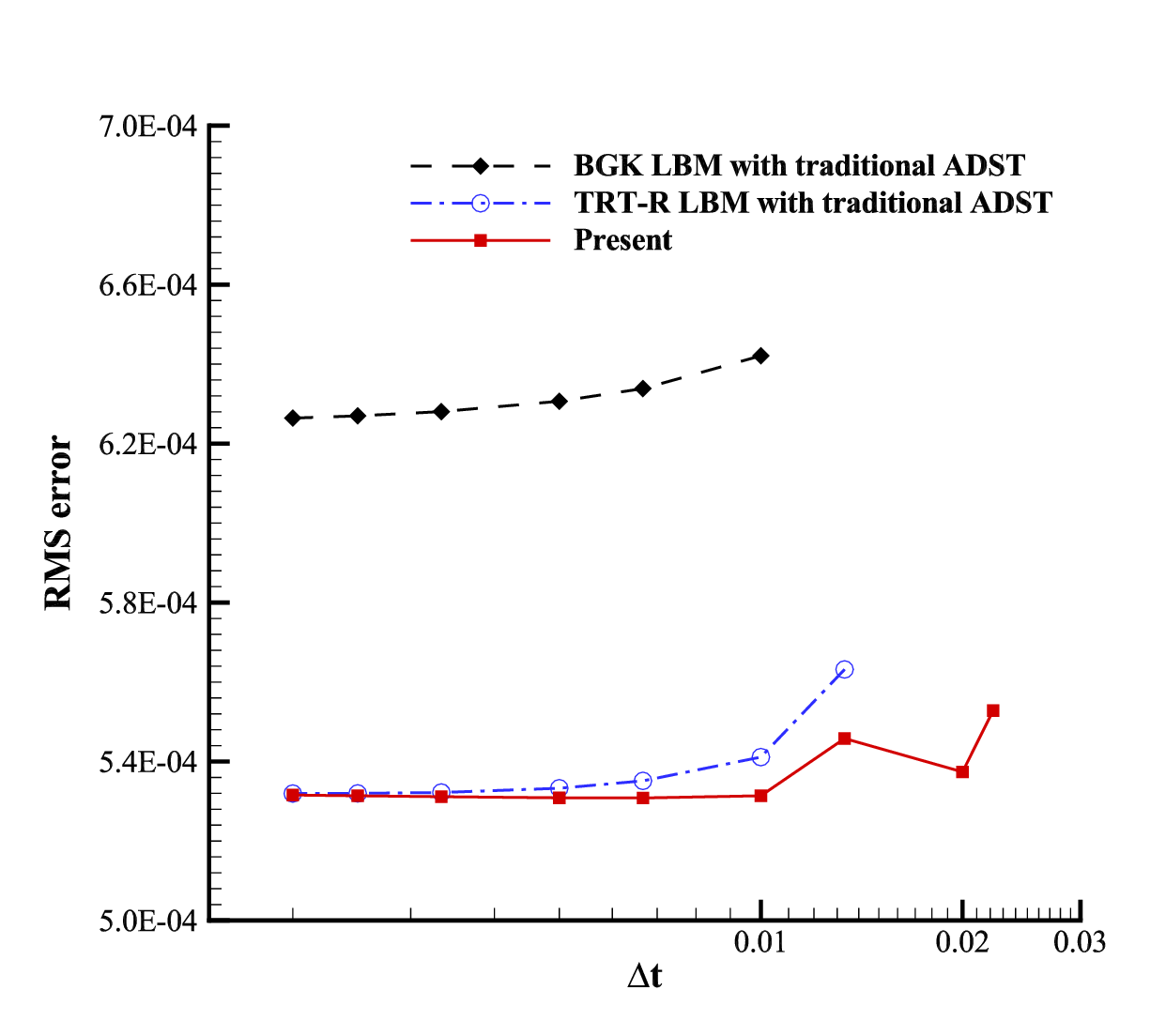}
 		\caption{$N^2=100^2$}
 		\label{fig1:sub1}
 	\end{subfigure}
 	\hfill
 	\begin{subfigure}[b]{0.45\linewidth}
 		\includegraphics[width=\linewidth]{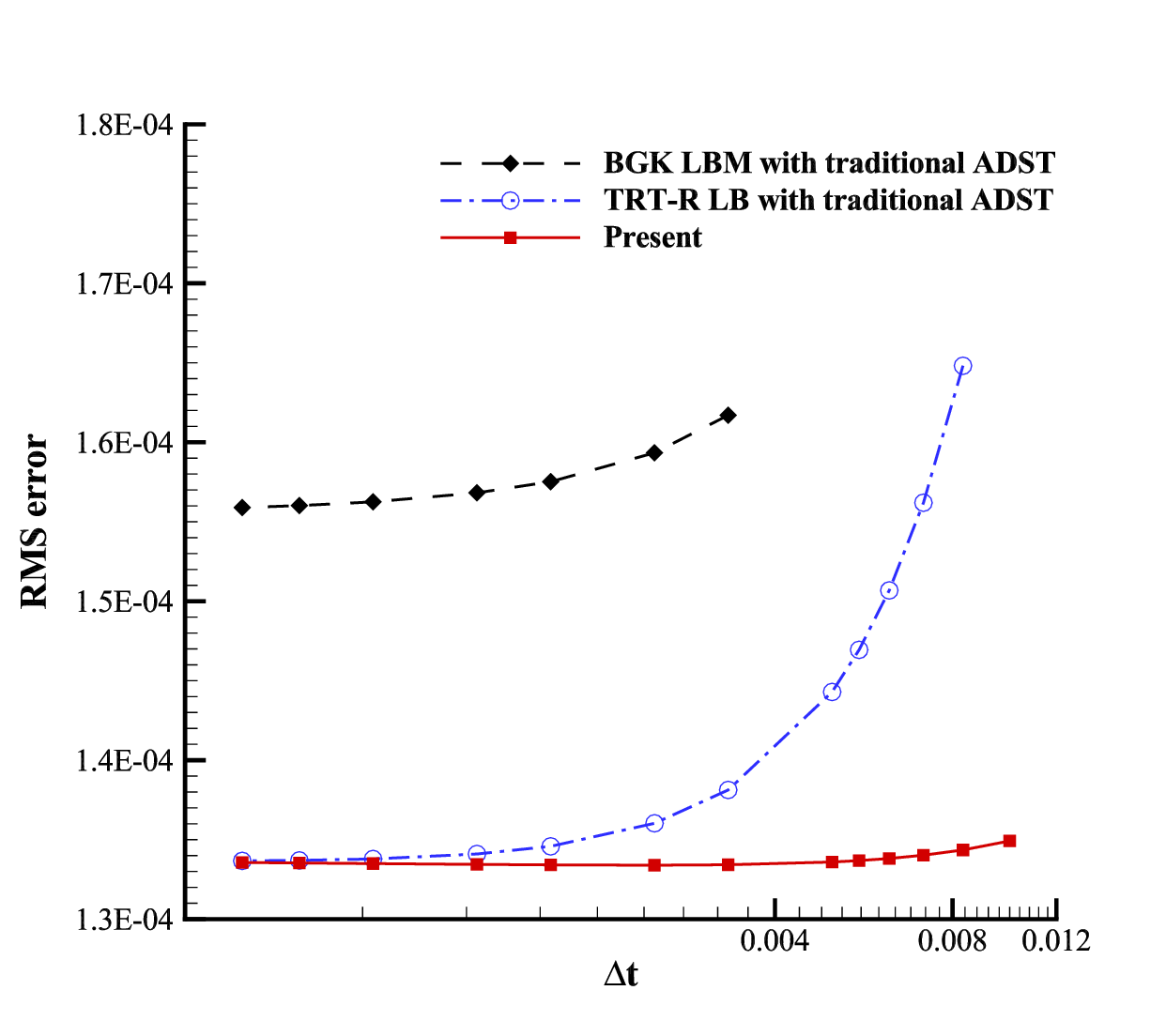}
 		\caption{$N^2=200^2$}
 		\label{fig1:sub2}
 	\end{subfigure}
 	\\
 	\begin{subfigure}[b]{0.45\linewidth}
 		\includegraphics[width=\linewidth]{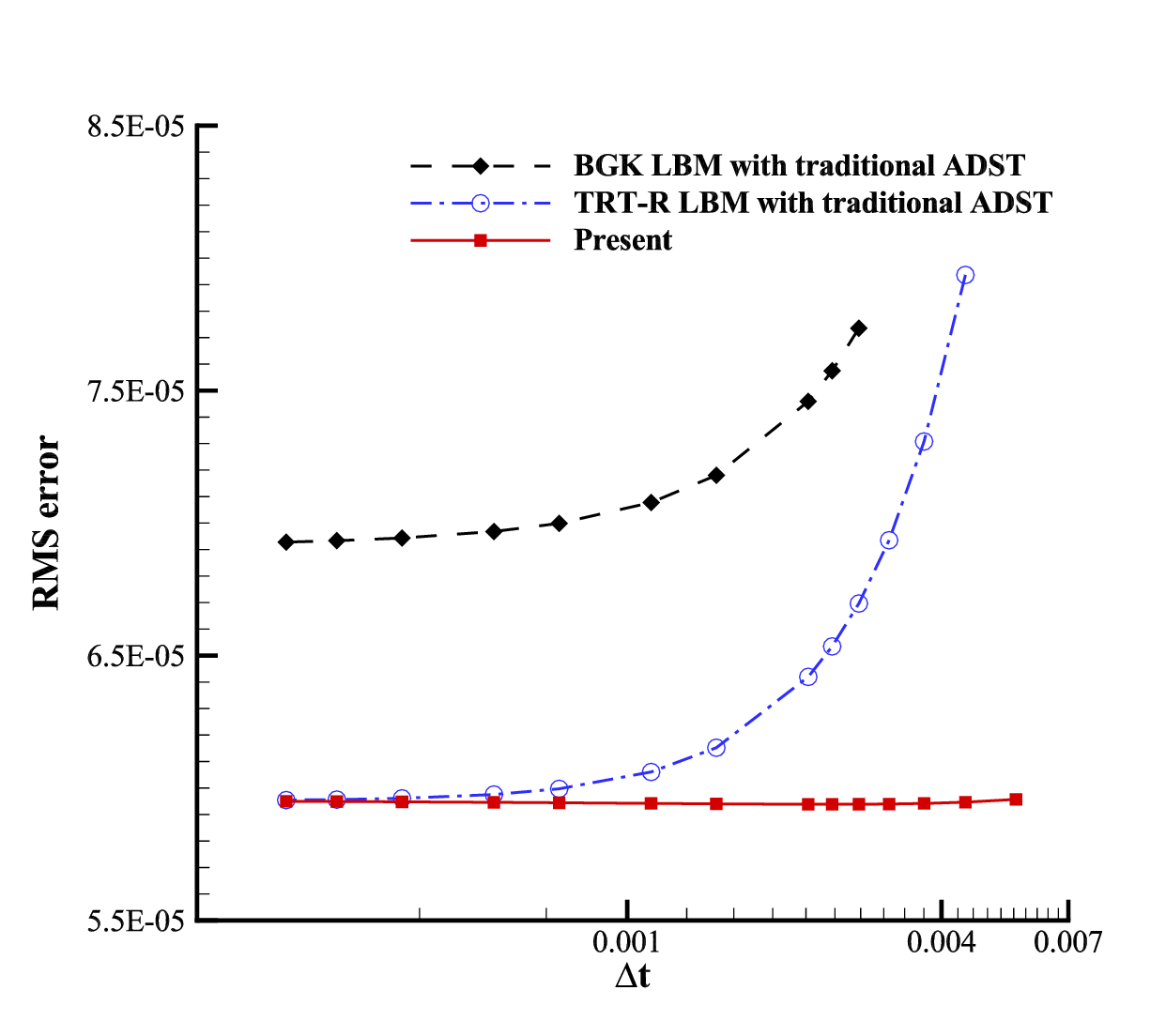}
 		\caption{$N^2=300^2$}
 		\label{fig1:sub3}
 	\end{subfigure}
 	\hfill
 	\begin{subfigure}[b]{0.45\linewidth}
 		\includegraphics[width=\linewidth]{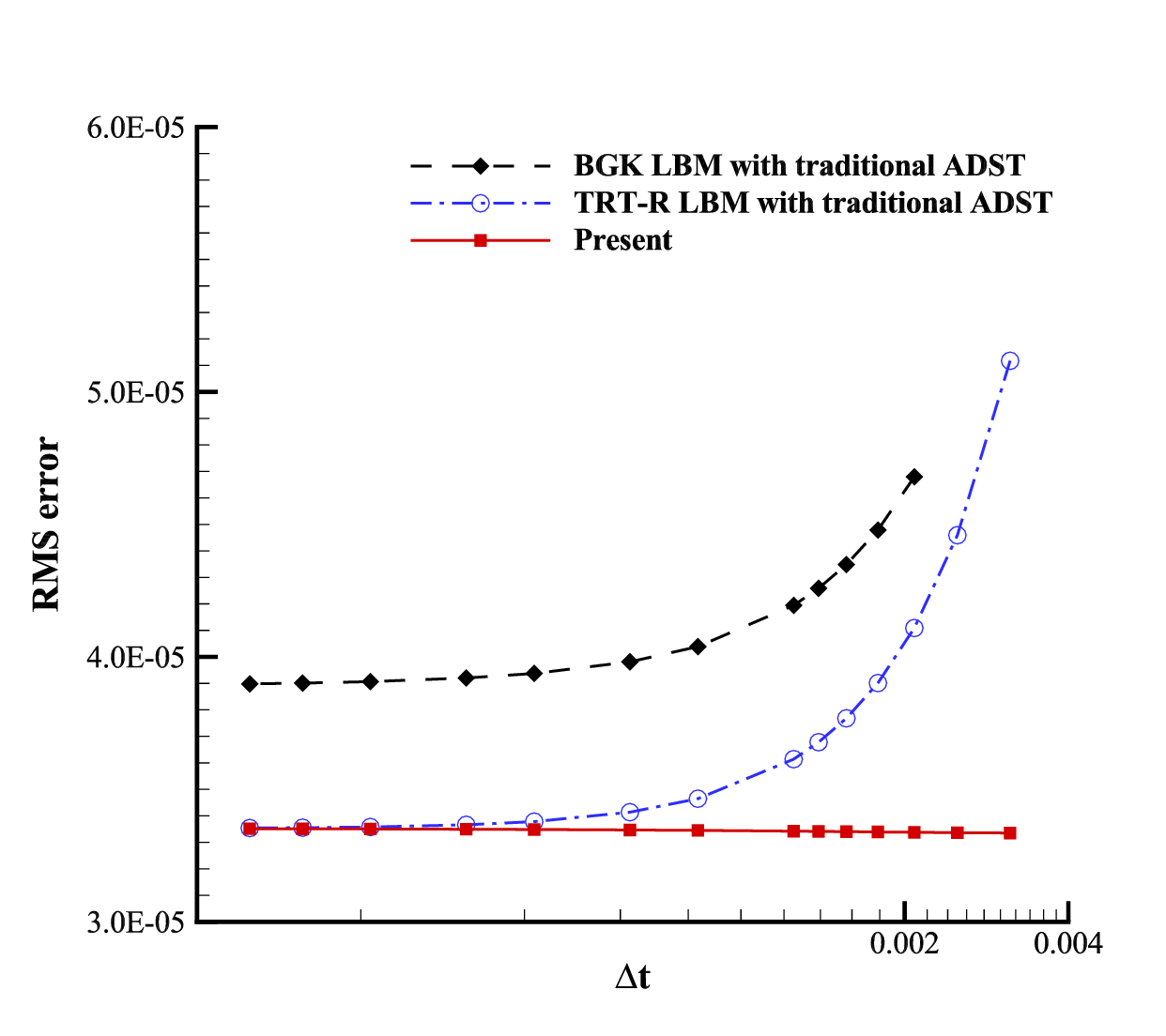}
 		\caption{$N^2=400^2$}
 		\label{fig1:sub4}
 	\end{subfigure}
 	\captionsetup{labelsep=period,labelfont=bf,labelformat=myformat}
 	\caption{Comparisons of the RMS errors $E_r$ of three models at different time steps with $\kappa=0.005$ for four kinds of grid resolutions. }
 	
 	\label{fig1}
 \end{figure}

\begin{figure}[H]
	\centering
	\begin{subfigure}[b]{0.45\linewidth}
		\includegraphics[width=\linewidth]{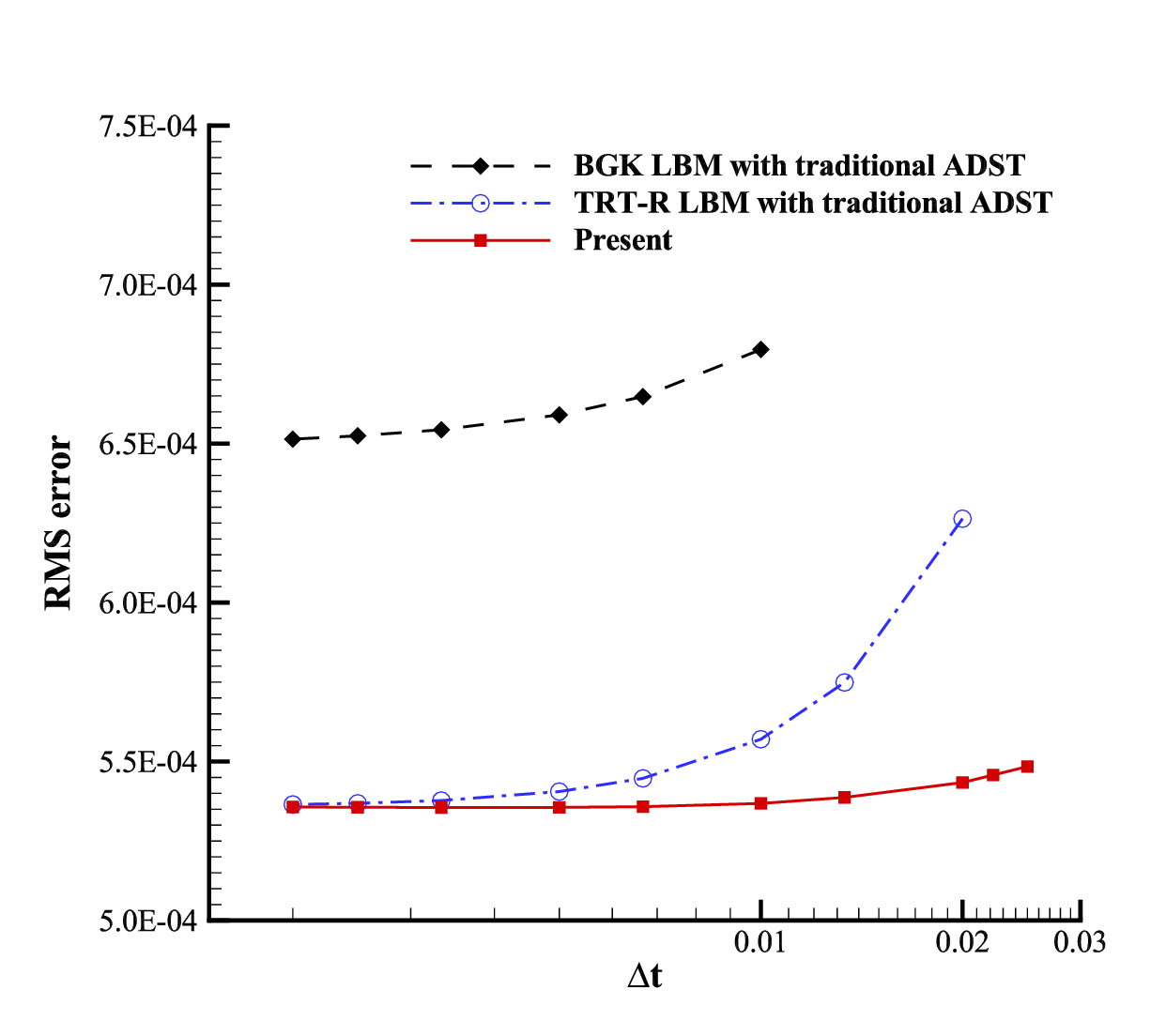}
		\caption{$N^2=100^2$}
		\label{fig2:sub1}
	\end{subfigure}
	\hfill
	\begin{subfigure}[b]{0.45\linewidth}
		\includegraphics[width=\linewidth]{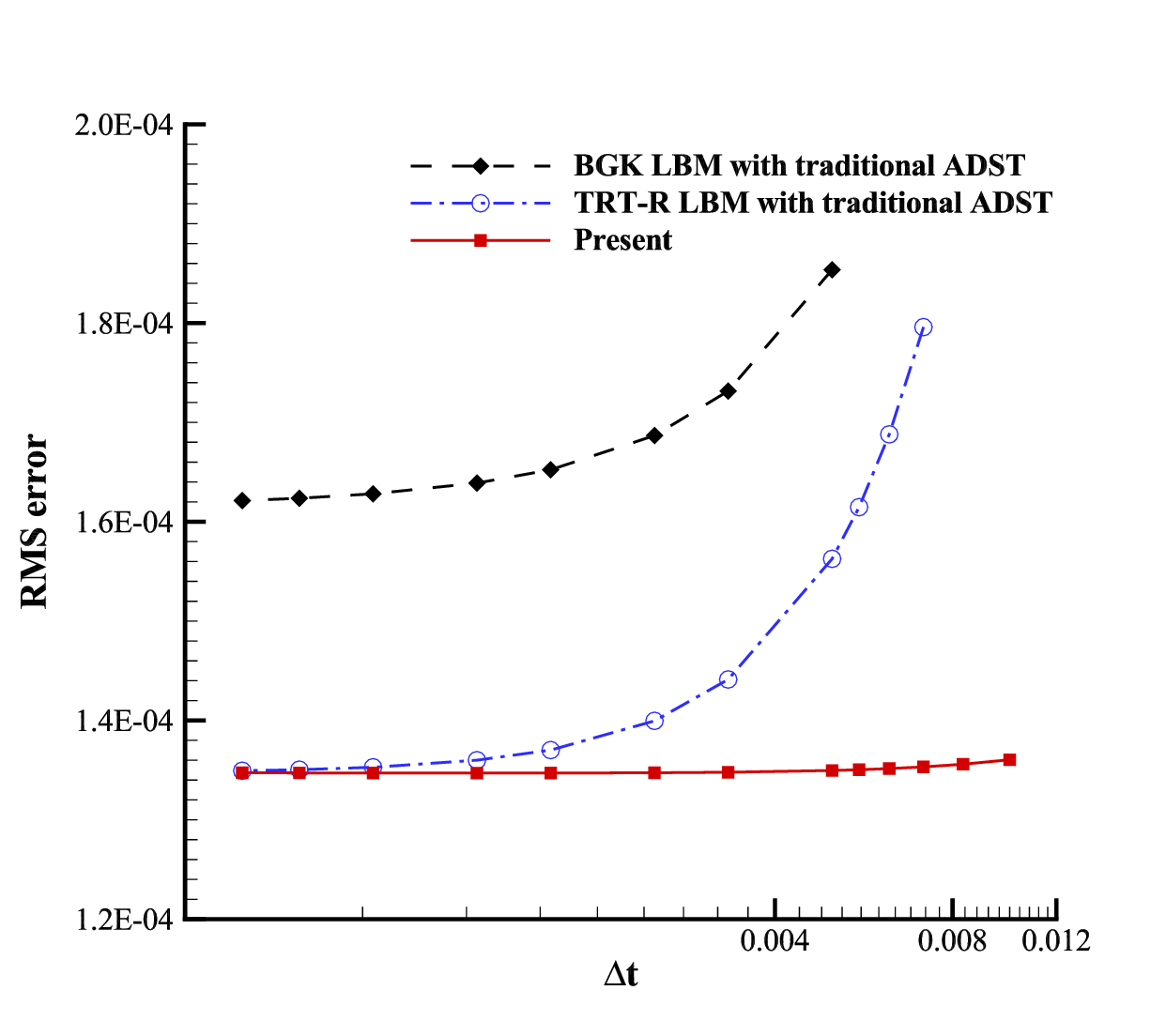}
		\caption{$N^2=200^2$}
		\label{fig2:sub2}
	\end{subfigure}
	\\
	\begin{subfigure}[b]{0.45\linewidth}
		\includegraphics[width=\linewidth]{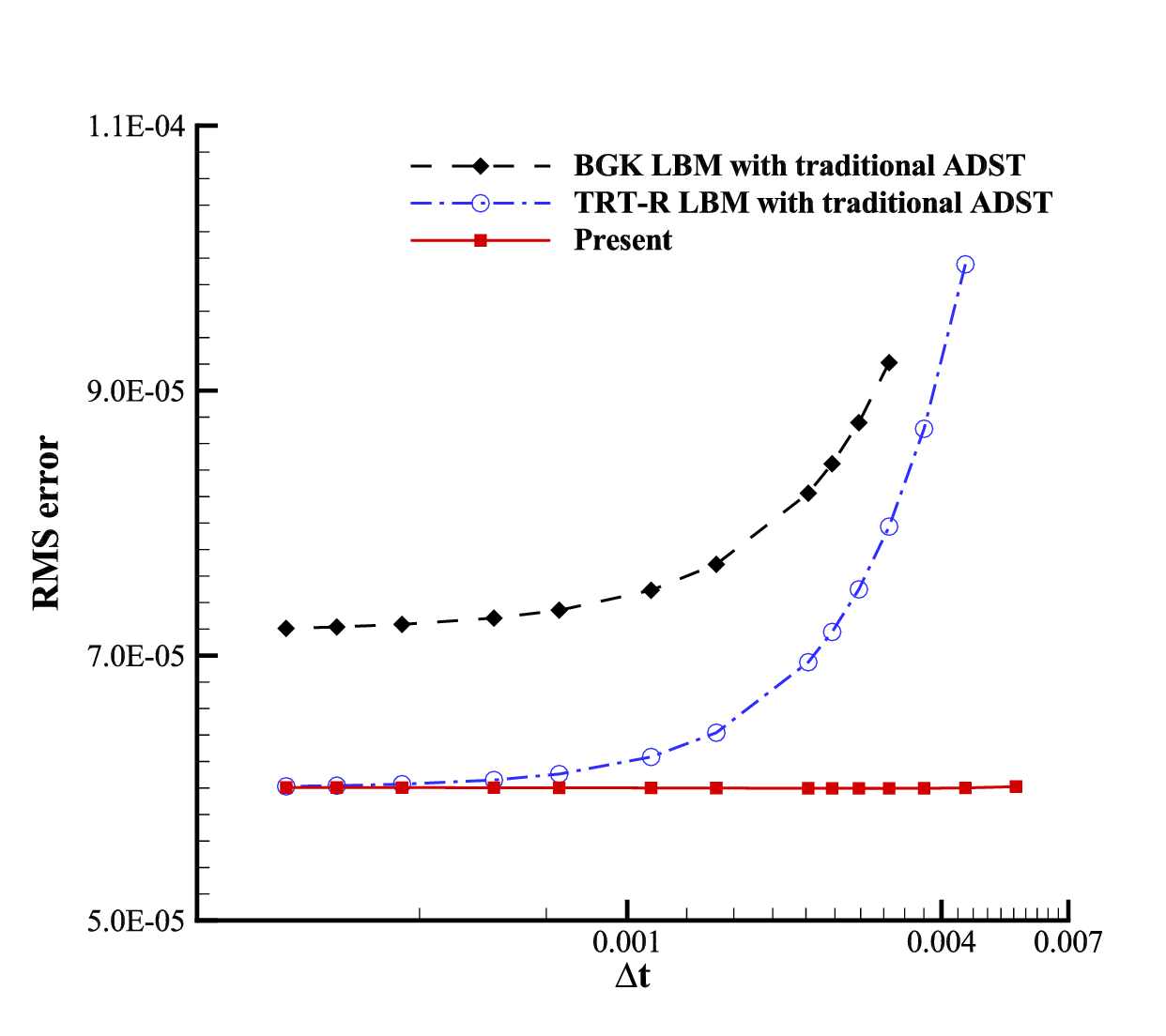}
		\caption{$N^2=300^2$}
		\label{fig2:sub3}
	\end{subfigure}
	\hfill
	\begin{subfigure}[b]{0.45\linewidth}
		\includegraphics[width=\linewidth]{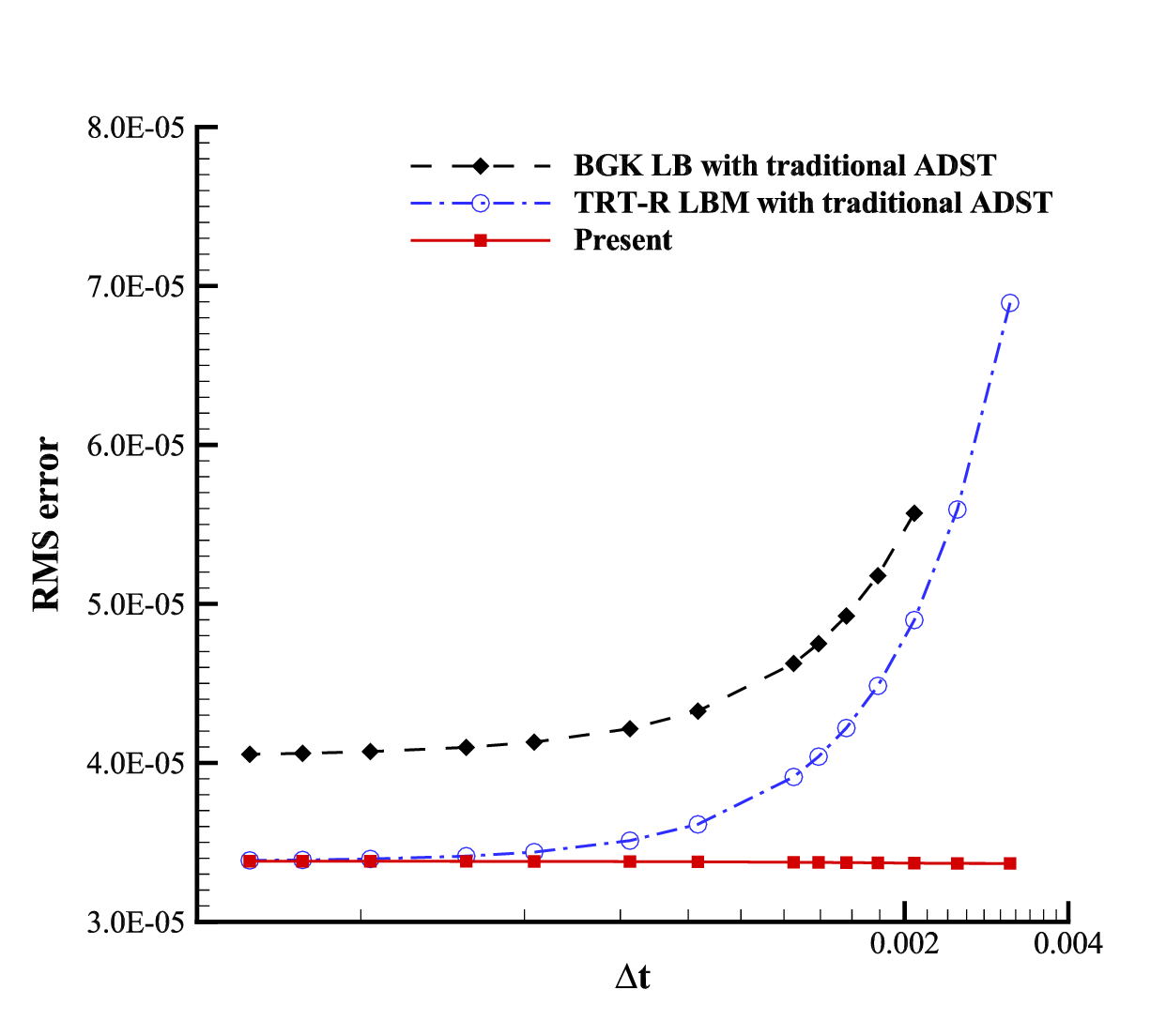}
		\caption{$N^2=400^2$}
		\label{fig2:sub4}
	\end{subfigure}
	\captionsetup{labelsep=period,labelfont=bf,labelformat=myformat}
	\caption{Comparisons of the RMS errors $E_r$ of three models at different time steps with $\kappa=0.01$ for four kinds of grid resolutions. }
	
	\label{fig2}
\end{figure}

\begin{figure}[H]
	\centering
	\begin{subfigure}[b]{0.45\linewidth}
		\includegraphics[width=\linewidth]{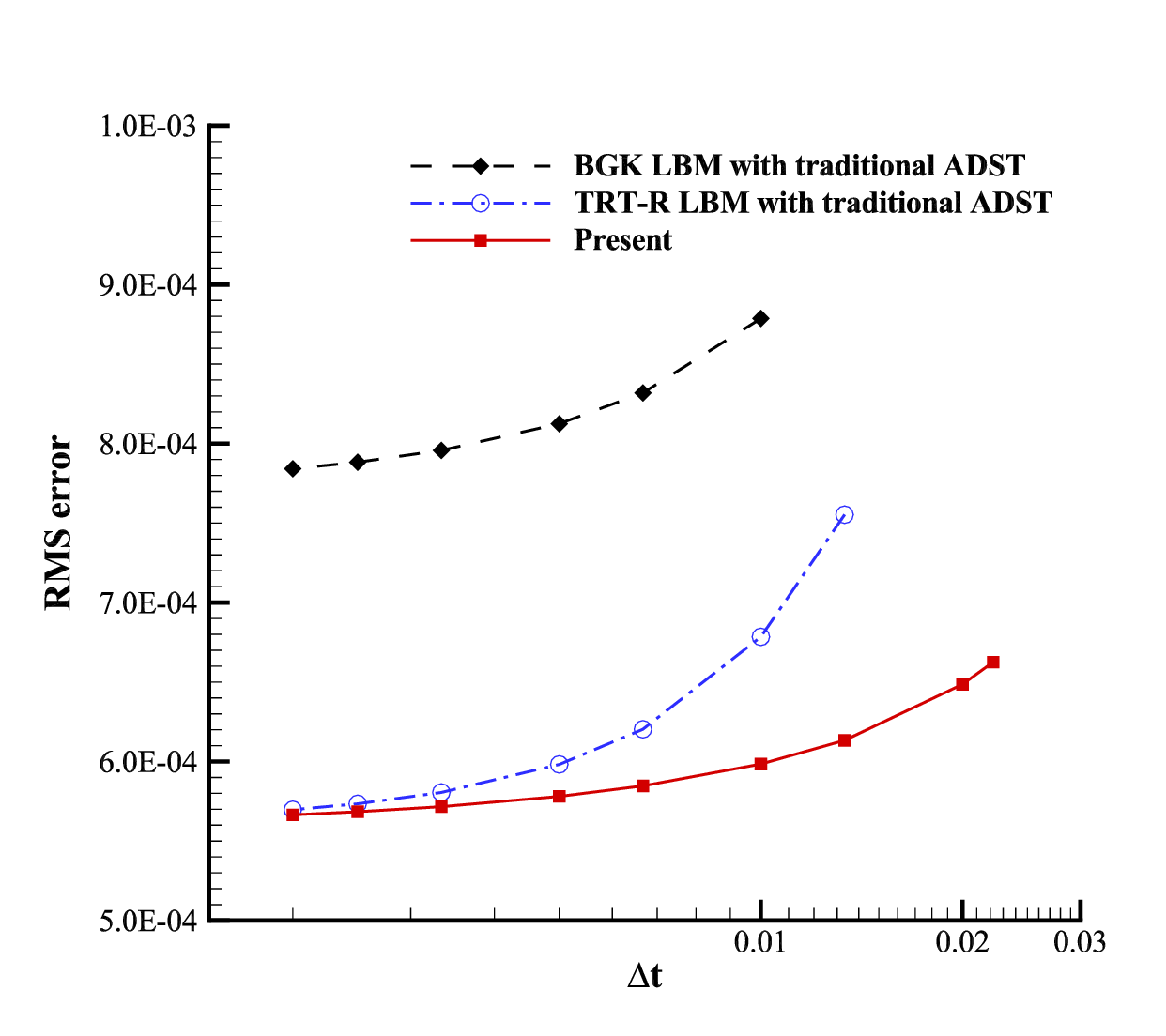}
		\caption{$N^2=100^2$}
		\label{fig3:sub1}
	\end{subfigure}
	\hfill
	\begin{subfigure}[b]{0.45\linewidth}
		\includegraphics[width=\linewidth]{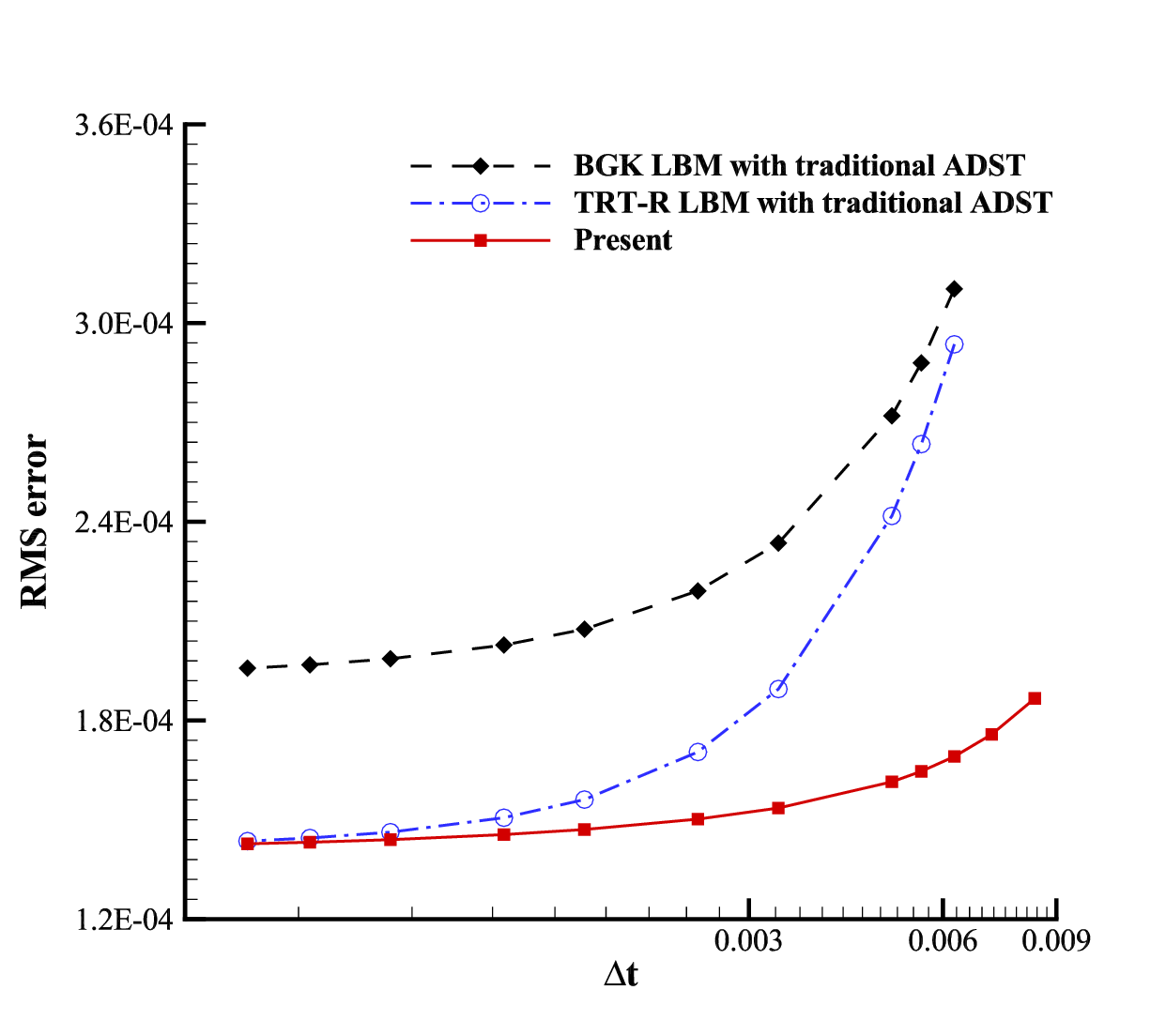}
		\caption{$N^2=200^2$}
		\label{fig3:sub2}
	\end{subfigure}
	\\
	\begin{subfigure}[b]{0.45\linewidth}
		\includegraphics[width=\linewidth]{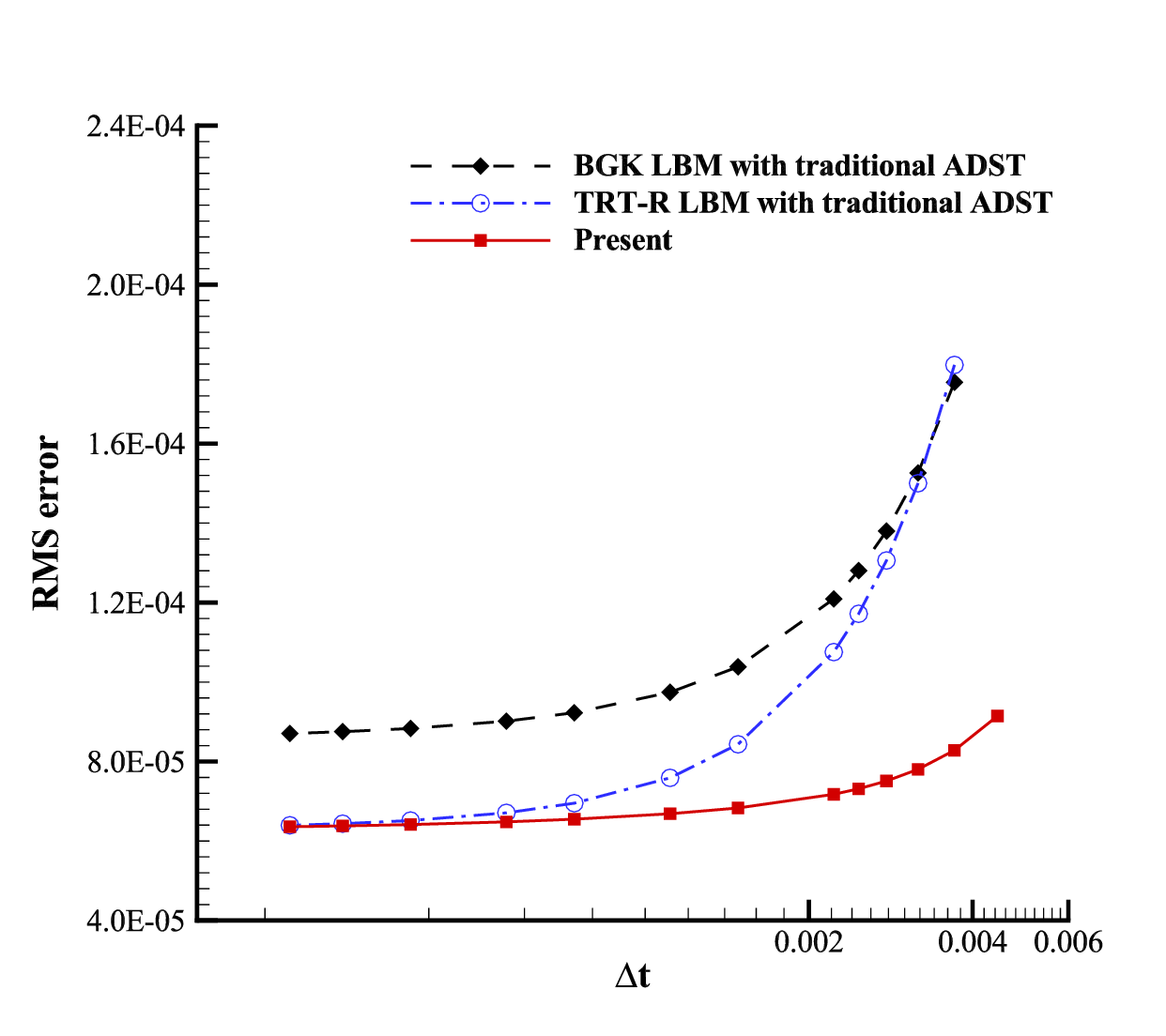}
		\caption{$N^2=300^2$}
		\label{fig3:sub3}
	\end{subfigure}
	\hfill
	\begin{subfigure}[b]{0.45\linewidth}
		\includegraphics[width=\linewidth]{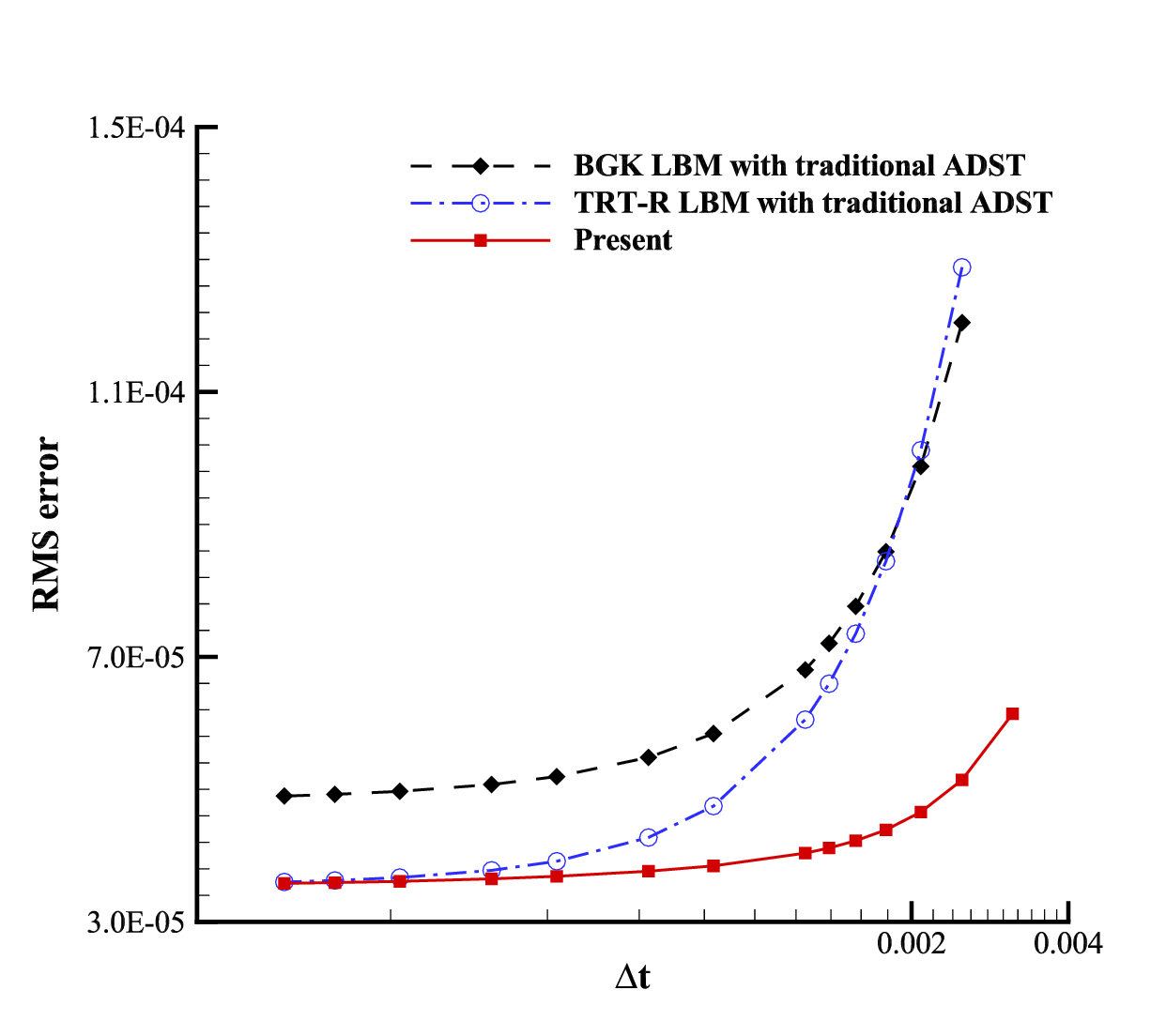}
		\caption{$N^2=400^2$}
		\label{fig3:sub4}
	\end{subfigure}
	\captionsetup{labelsep=period,labelfont=bf,labelformat=myformat}
	\caption{Comparisons of the RMS errors $E_r$ of three models at different time steps with $\kappa=0.05$ for four kinds of grid resolutions. }
	
	\label{fig3}
\end{figure}

Now we turn to the accuracy of these models. In simulations, the following parameters: $N^2=100^2, 200^2, 300^2, 400^2$, $\kappa=0.05,0.01,0.005$, and several ${\Delta t}$ were considered.  Fourteen different time steps were used for each grid size as shown in Table \ref{table_dt}. The RMS error of each example at $t=1$ was recorded if the program was stable. The results for $\kappa=0.005$, $0.01$, and $0.05$ are shown in \hyperref[fig1]{Fig.~\ref{fig1}}, \hyperref[fig2]{Fig.~\ref{fig2}}, and \hyperref[fig3]{Fig.~\ref{fig3}}, respectively. As we observed, the previous models, i.e., the BGK LBM with traditional ADST\cite{malaspinas2010lattice} and the TRT-R LBM with traditional ADST\cite{zhao2020block}, tolerate increasing numerical errors as the time step increases, leading to blow-up at large time steps. But the present model always exhibits the lowest error and the best stability.  

\begin{table}
	\captionsetup{
		labelsep=newline, 
		justification=raggedright, 
		singlelinecheck=false, 
		labelformat=Table_fm 
	}
	\caption{Settings of $\Delta t$ for different grid sizes $N^2$.}
	\label{table_dt}
	\begin{tabular}{lllll}
		\toprule
		Grid sizes & $N^2=100^2$ & $N^2=200^2$ & $N^2=300^2$ &$N^2=400^2$  \\
		\midrule
		& 1/20 & 1/80 	& 1/180 & 1/320 \\
		& 1/25 & 1/100 	& 1/225 & 1/400 \\
		& 1/30 & 1/120 	& 1/270 & 1/480 \\
		& 1/35 & 1/140 	& 1/315 & 1/560 \\
		& 1/40 & 1/160 	& 1/360 & 1/640 \\
		& 1/45 & 1/180 	& 1/405 & 1/720 \\
		$\Delta t$ & 1/50 & 1/200 	& 1/450 & 1/800 \\
		& 1/75 & 1/300 	& 1/675 & 1/1200 \\
		& 1/100 & 1/400 & 1/900 & 1/1600 \\
		& 1/150 & 1/600 & 1/1350 & 1/2400 \\
		& 1/200 & 1/800 & 1/1800 & 1/3200 \\
		& 1/300 & 1/1200 & 1/2700 & 1/4800 \\
		& 1/400 & 1/1600 & 1/3600 & 1/6400 \\
		& 1/500 & 1/2000 & 1/4500 & 1/8000 \\
		\bottomrule
	\end{tabular}
\end{table}

\subsection{Convergence rate}
Furthermore, the convergence rate of the present model is investigated. Several grid sizes $N^2=20^2, 40^2, 60^2, 80^2, 100^2$ were considered here. We set $\frac{(\Delta x)^2}{\Delta t}=0.16\pi^2$ to keep the relaxation time $\tau_{p,1}$ constant for a given diffusion coefficient $\kappa$. Results for diffusion coefficients (e.g. $\kappa=0.005$, $\kappa=0.01$, and $\kappa=0.05$) are shown in \hyperref[table1]{Table~\ref{table1}}. It can be observed that the BGK model undergoes blow-up when $N<100$, whereas the other two models exhibit good convergence. Specifically, the present model has a slightly higher convergence rate than the traditional model of Zhao et al.\cite{zhao2020block}. 

\begin{table}[!htbp]
	\captionsetup{
		labelsep=newline, 
		justification=raggedright, 
		singlelinecheck=false, 
		labelformat=Table_fm 
	}
	\caption{Comparision of the convergence rates of different models with $\kappa=0.01$ at $t=1.0$. (Model 1: TRT-R LBM with traditional ASDT\cite{zhao2020block}, Model 2: BGK LBM with traditional ASDT\cite{malaspinas2010lattice}. ``$\circ$'' and ``$-$'' indicate unstable and incalculable order, repectively.)}
	\label{table1}
	\begin{tabular}{*{8}{l}} 
		\toprule
		\multirow{2}*{Grid sizes} & \multicolumn{2}{l}{Present} & \multicolumn{2}{l}{Model 1} & \multicolumn{2}{l}{Model 2} & \\
		\cmidrule(lr){2-3} \cmidrule(lr){4-5} \cmidrule(lr){6-8}
		& RMS & Order & RMS & Order & RMS & Order \\
		\midrule

$\kappa=0.005, N=20$  & 1.196E-02 & $-$  & 1.182E-02 & $-$  & $\circ$      & $-$\\
$\kappa=0.005, N=40$  & 3.471E-03 & 1.676 & 3.767E-03 & 1.650 & $\circ$      & $-$\\
$\kappa=0.005, N=60$  & 1.442E-03 & 2.350 & 1.464E-03 & 2.331 & $\circ$      & $-$\\
$\kappa=0.005, N=80$  & 8.187E-04 & 1.969 & 8.330E-04 & 1.960 & $\circ$ 	 & $-$\\
$\kappa=0.005, N=100$ & 5.314E-04 & 1.937 & 5.411E-04 & 1.933 & 6.420E-04 & $-$\\

$\kappa=0.01, N=20$   & 1.110E-02 & $-$  & 1.118E-02 & $-$  & $\circ$      & $-$\\
$\kappa=0.01, N=40$   & 3.351E-03 & 1.729 & 3.442E-03 & 1.700 & $\circ$      & $-$\\
$\kappa=0.01, N=60$   & 1.473E-03 & 2.028 & 1.524E-03 & 2.010 & $\circ$      & $-$\\
$\kappa=0.01, N=80$   & 8.377E-04 & 1.961 & 8.683E-04 & 1.955 & $\circ$      & $-$\\
$\kappa=0.01, N=100$  & 5.368E-04 & 1.994 & 5.570E-04 & 1.990 & 6.796E-04 & $-$\\

$\kappa=0.05, N=20$   & 1.288E-02 & $-$  & 1.420E-02 & $-$  & $\circ$      & $-$\\
$\kappa=0.05, N=40$   & 3.679E-03 & 1.807 & 4.146E-03 & 1.776 & $\circ$      & $-$\\
$\kappa=0.05, N=60$   & 1.653E-03 & 1.974 & 1.870E-03 & 1.964 & $\circ$      & $-$\\
$\kappa=0.05, N=80$   & 9.332E-04 & 1.987 & 1.057E-03 & 1.982 & $\circ$      & $-$\\
$\kappa=0.05, N=100$  & 5.984E-04 & 1.991 & 6.785E-04 & 1.988 & 8.787E-04 & $-$\\
		\bottomrule
	\end{tabular}
	
\end{table}
\subsection{The optimal free relaxation parameter to eliminate the numerical slip of bounce-back scheme}
As mentioned ealier, the TRT CO with a specific free relaxation parameter can effectively eliminate the numerical slip of the bounce-back boundary scheme. The work of Zhao et al.\cite{zhao2020block} has shown that their B-TriRT model can indeed successfully eliminate the numerical slip, and although the TRT-R collision operator in this paper can be regarded as a special case of the B-TriRT collision operator, a concrete numerical verification is still necessary.

\begin{figure}[H]
	\centering
	\includegraphics[width=0.75\textwidth]{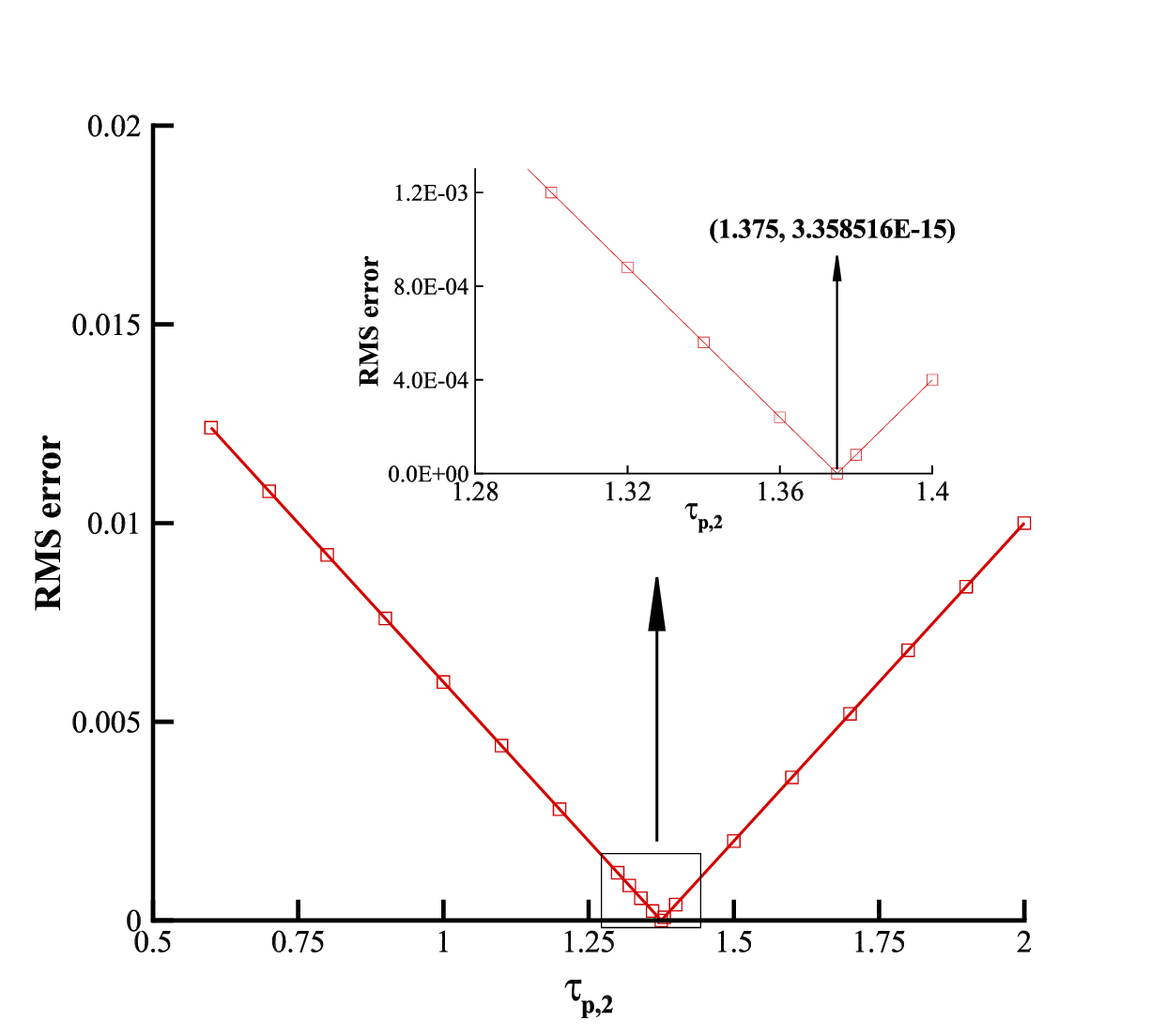}
	\captionsetup{labelsep=period,labelfont=bf,labelformat=myformat}
	\caption{The RMS errors at different free relaxation parameters $\tau_{p,2}$ with $\tau_{p,1}=0.8$ and theoretical $\tau_{p,2}^t=1.375$. }
	
	\label{fig4}
\end{figure}

Using the parameter settings and convergence criteria corresponding to Fig.~3(b) in the paper of Zhao et al.\cite{zhao2020block}, we can simulate the unidirectional steady-state diffusion problem. Note that the Dirichlet boundary conditions need to be executed by the halfway bounce-back scheme on the top and bottom rest boundaries. According to the analysis of Zhao et al.\cite{zhao2020block}, the relaxation parameters in this paper should satisfy the following relation to eliminate the numerical slip: 
\begin{equation}\label{eq_tau}
	\tau_{p,2}=\frac{3(1-4\tau_{p,1})}{8(1-2\tau_{p,1})}.
\end{equation}
Here, we set $\tau_{p,1}=0.8$, so the theoretical ``free'' relaxation parameter is $\tau_{p,2}^t=1.375$ according to the above equation. We recorded the RMS errors at different free relaxation parameters $\tau_{p,2}$ and showed them in Fig.~\ref{fig4}. The results show that, for the halfway bounce-back Dirichlet boundary conditions, the derived TRT-R CO can indeed eliminate the numerical slip with the relaxation parameter described by Eq.~(\ref{eq_tau}).
\section{Conclusion}\label{section4}
In this paper, there are three main improvements in terms of model construction. Fistly, a new first-order space-derivative auxiliary term was proposed instead of the traditional time-derivative term to derive the correct CDE without the help of NSEs. Secondly, the first-order discrete-velocity form of discrete source term was adopted to improve the accuracy of the source term. Finally, the TRT-R CO was derived by the non-equilibrium high-order Hermite expansion. Based on these improvements, we proposed the TRT-R LB model for the CDE with variable coeffients, which is often found in the solution of the CDE within a time step coupled with the equation controlling the convection velocity. The convection and diffusion of Gaussian pluse in a rigid-body rotating flow was simulated by the proposed model and two other popular models, i.e., the TRT-R model with traditional ADST\cite{zhao2020block} and the BGK model with traditional ADST\cite{malaspinas2010lattice}. Results show that our model can effectively solve the CDE systems with non-uniform convection velocity, and has better accuracy, stability and convergence than the previous models, especially for a large time step. This also hints at the potential of the proposed space-derivative auxiliary term in convection-dominant problems.

It should be noted that, although the B-TriRT CO can be transformed to the TRT-R CO by choosing specific relaxation parameters, we presented an alternative way of deriving the TRT-R CO. According to this idea, one can construct more types of regularized collision operators by using non-equilibrium higher-order Hermite expansion and adjusting the free relaxation time parameters. Besides, it can be found that the TRT-R CO has better stability and accuracy than the BGK CO and can be used to eliminate the numerical slip of BB scheme. However, the optimal free relaxation parameter under different circumstances of the TRT-R CO is beyond the scope of this paper and will be clarified in future work. 

\section*{Acknowledgments}
This work is financially supported by the National Natural Science Foundation of China (Grant
Nos. 12101527, 12271464 and 11971414), the Science and Technology Innovation Program of Hunan Province (Program No. 2021RC2096), Project of Scientific Research Fund of Hunan Provincial Science and Technology Department (Grant No. 21B0159) and the Natural Science Foundation for Distinguished Young Scholars of Hunan Province (Grant No. 2023JJ10038). 

\section*{Declaration of competing interest}
The authors have no conflicts to disclose.

\bibliographystyle{elsarticle-num} 
\bibliography{mybibfile}

%
%
%
%
\end{document}